\documentclass[11pt,reqno]{amsart}

\usepackage[T1]{fontenc}
\usepackage{graphicx}
\usepackage{color}
\definecolor{MyLinkColor}{rgb}{0,0,0.4}

\newcommand{\gr}{\mathop{\rm graph}\nolimits}

\newcommand{\0}{\Omega}

\newcommand{\p}{\partial}

\newcommand{\R}{\mathbb{R}}

\newtheorem{thm}{Theorem}[section]

\theoremstyle{remark} 
\newtheorem{rem}[thm]{Remark}

\numberwithin{equation}{section}

\title[On particle trajectories in linear deep-water waves]{On particle trajectories in linear deep-water waves}


\usepackage[colorlinks=true,linkcolor=MyLinkColor,citecolor=MyLinkColor]{hyperref}
\author[A.-V. Matioc]{Anca-Voichita Matioc}
\address{Institut f{\"u}r Angewandte Mathematik, Leibniz Universit{\"a}t Hannover, Welfengarten~1, 30167 Hannover, Germany. }
\email{matioca@ifam.uni-hannover.de}

\usepackage[colorlinks=true,linkcolor=MyLinkColor,citecolor=MyLinkColor]{hyperref} %dvips

\begin{document}

\begin{abstract}
We determine  the phase portrait of a Hamiltonian system of equations describing
the motion of the particles in linear deep-water waves.  
The particles experience in each period a forward drift which decreases with greater depth.
\end{abstract}

\maketitle

\section{Introduction}
The motion of water particles under the waves which advance across the water is a classical problem in this field.
Watching the sea it is oft possible to trace a wave as it propagates on the water's surface,
 but what one observes traveling across the sea is not the water but a wave pattern.
Although the wave travels from one place to another, the substance through which it travels moves very
little. As the wave advances across the water and can be followed for a long way, a typical water particle moves slightly
up and down, forward and backward as the wave passes it.
 If an object hovered in the water, like a water particle, its motion
would be synchronized with that of a floating object lying on the water's surface, with its orbit diminishing with the
distance from the surface. As waves generated by wind in an area move towards a region where the wind has ceased,
we observe swell-long crested two-dimensional waves approaching a smooth sinusoidal shape and moving over long
distances.  
%Swell enters the category of deep-water waves (waves in water of depth greater than half of the wavelength). 
Deep-water waves are modelled mathematically as periodic two-dimensional waves in water of infinite
depth. The motion of the water particle in the fluid below swell is of great interest.
The classical description of these
particle paths is obtained within the framework of linear water wave theory \cite{Ac,Cr,DE, RJ,LA,LI, So, St}: all water particles
trace a circular orbit, the diameter of which decreases with depth so that the orbital motion practically ceases at depth
equal to one-half the wavelength. These features have important practical consequences. For example, a submarine
at a depth below half a wavelength would hardly notice the motion of the surface wave, for this reason submarines
dive during storms in the open sea.

The only known explicit solution with a non-flat free surface of the governing
equations for gravity water waves is Gerstner's wave \cite{GE}: a deep-water wave solution for which all particle paths
are circles of diameters decreasing with the distance from the free surface (see the discussion in \cite{C, C1}). Due to the
mathematical intractability of the governing equations for water waves, for irrotational water waves (Gerstner's wave
has a peculiar nonvanishing vorticity) the classical approach \cite{Cr, HD2, RJ,LA,MT,So,St} relies on analyzing the particle
motion after linearization of the governing equations. However, even within the linear water wave theory, the ordinary
differential equations system describing the motion of the particles is nevertheless nonlinear and explicit solutions of
this system are not available, but qualitative features of the underlying flow are known \cite{CE,CS, CS1, CSat, Hur, AMat, Mati, TO} and could possibly lead to a confirmation of the features we prove here, within the confines of linear water waves (see \cite{C2}).

In this paper we show  that for linear deep-water wave no particle trajectory is actually closed, unless the free surface is flat. 
Each trajectory involves over a
period a backward/forward movement of the particle, and the path is an elliptical
arc  with a forward drift. Support for
this conclusion is given by the analysis of the average flow of energy within linear water wave theory (see \cite{RJ}): due to
the passage of a periodic surface wave, the water particles in the fluid experience on average a net displacement in the
direction in which the waves are propagating, termed Stokes drift \cite{Sts} (see the discussion in \cite{Y}).
This result was first obtained  in \cite{CEV} for gravity waves and \cite{HD1} for capillary and capillary-gravity water waves.
This results have been  improved in \cite{C2, CS, DHH} in the case of irrotational waves by describing the exact geometry  of the actual particle paths. 
For gravity deep-water waves we reconfirm herein the results of \cite{CEV}, but our  analysis is more precise and uses  only elementary analysis methods.
The work is structured as follows: we recall first  the governing equations for water waves and 
give the exact solutions of the linearised problem.
Our major contribution, in Section 3, is a precise  analysis of the phase portrait of a Hamiltonian system associated to the 
particle motion.
This is the key point in determining the trajectories of the water particles.

\section{Preliminaries}
In this section we recall the governing equations for the propagation of two-dimensional gravity deep-water waves and we present
their linearization (for a more detailed discussion we refer to \cite{RJ}).

\subsection{The governing equations}
We consider a two-dimensional inviscid incompressible fluid in a constant gravitational field.
To describe the waves propagation 
we consider a cross section of the flow that is perpendicular to the crest line
with Cartesian coordinates $(x, y)$, the $y$-axis pointing vertically upwards and the
$x$-axis being the direction of wave propagation (see Figure 1). Let $(u(t, x, y), v(t, x, y))$ be the velocity field and  the equation of the free surface with mean water level zero is $y =\eta(t, x) $ with  $\int_\R \eta(t,x)=0$.

Constant density (homogeneity) is a physically reasonable assumption for gravity waves \cite{RJ}, and it implies the equation of mass conservation
\begin{equation}\label{mc}
u_x+u_y=0.
\end{equation} 
Under the assumption of inviscid flow the equation of motion is Euler's equation
\begin{equation}\label{Euler}
\left\{
\begin{array}{rllllll}
u_t+uu_x+vu_y&=&-P_x \\[2ex]
v_t+uv_x+vv_y&=& -P_y-g,
\end{array}
\right.
\end{equation}
where $P(t,x,y)$ denotes the pressure and $g$ is the gravitational constant of 
acceleration. The free surface decouples the motion of the water from that of the air,
a fact that is expressed in \cite{RJ} by the dynamic boundary condition
\begin{equation}\label{Dbc}
P=P_0 \quad \quad \text{on} \quad y=\eta(t,x),
\end{equation}
if we neglect surface tension, where $P_0$ is the (constant) atmospheric pressure.
Since the same particles always form the free surface, we also have the kinematic
boundary condition
\begin{equation}\label{kbc}
v=\eta_t+u\eta_x \quad \quad \text{on} \quad y=\eta(t,x).
\end{equation}
The boundary condition at the bottom
\begin{equation}\label{boc}
(u,v) \to (0,0) \, \, \text{as} \, y\to -\infty \, \, \text{uniformly for}\, \, x\in \R,
\end{equation}
expresses  the fact that at great depths there is practically no motion. Summarising, the governing equations for water
waves are encompassed by the nonlinear free-boundary problem 
\begin{equation}\label{fbpb}
\left\{
\begin{array}{rllllllll}
u_x+u_y &=& 0, \\[2ex]
u_t+uu_x+vu_y&=&-P_x, \\[2ex]
v_t+uv_x+vv_y&=& -P_y-g, \\[2ex]
P&=&P_0 &\quad  \text{on} \quad y=\eta(t,x),\\[2ex]
v&=& \eta_t+u\eta_x &\quad  \text{on} \quad y=h_0+\eta(t,x),\\[2ex]
(u,v) &\to& (0,0)  &\quad  \text{as} \quad y \to -\infty\, \text{uniformly for}\, \, x\in \R.
\end{array}
\right.
\end{equation}

\begin{figure}\label{F:Water} 
$$\includegraphics[width=11cm]{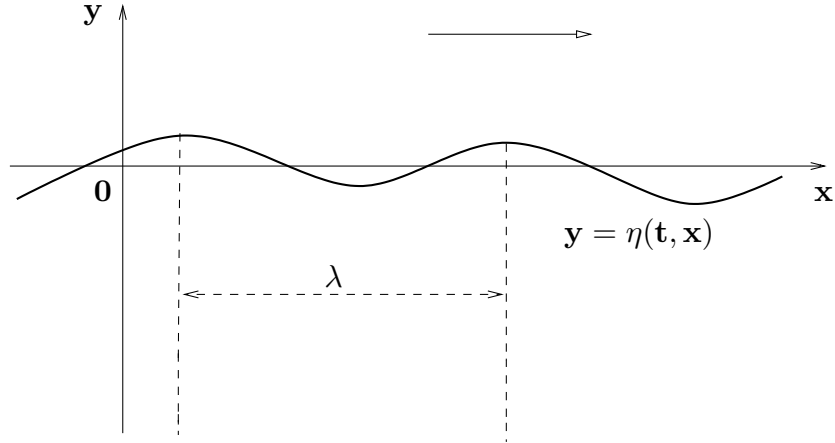}$$
\caption{A deep-water wave with wavelength $\lambda$}
\end{figure}

An important category of flows are those of zero vorticity (irrotational flows),
characterized by the additional assumption
\begin{equation}\label{ad}
u_y=v_x.
\end{equation}
The vorticity of a flow, $\omega = u_y - v_x$, measures the local spin or rotation of a
fluid element so that in irrotational flows the local whirl is completely absent.
Relation \eqref{ad} ensures the existence of a velocity potential $ \phi(t, x, y)$ defined up to
a constant via
\[
\phi_x=u, \quad \quad \phi_y=v.
\]
In view of \eqref{mc}, $\phi$ is a harmonic function, i.e. $(\partial_x^2+\partial_y^2)\phi=0$
so that the 
methods of complex analysis become available for the study of irrotational flows.

Concerning the physical relevance of irrotational water flows, field evidence 
indicates that for waves entering a region of still water the assumption of irrotational
flow is realistic \cite{LI}. Moreover, as a consequence of Kelvin's circulation theorem
\cite{Ac}, a water flow that is irrotational initially has to be irrotational at all later times.
It is thus reasonable to consider that water motions starting from rest will remain
irrotational at later times. Nonzero vorticity is the hallmark of the interaction between waves and non-uniform currents cf. \cite{CS1} and \cite{CV, IK}.
The most important example of rotational flows are tidal currents cf. the discussion in \cite{C3, TP}, but the assumption of constant vorticity can not be accomodated within the framework of deep-water flows. However, we consider herein irrotational gravity waves. 

\subsection{Linear water waves}
The problem \eqref{fbpb} is nondimensionalised using a typical wavelength $\lambda$ and a typical amplitude of the wave $\varepsilon$.
We define the set of nondimensional variables 
\begin{align*}
&x\mapsto \lambda x, \quad y\mapsto y, \quad t\mapsto \frac{\lambda}{\sqrt{g}}t,\quad u\mapsto u\sqrt{g}, \quad v\mapsto v\frac{\sqrt{g}}{\lambda}, \quad \eta\mapsto \varepsilon\eta,
\end{align*}
combined with the scaling 
\[
p \mapsto  \varepsilon p, \quad (u, v) \mapsto \varepsilon (u, v).
\]
 Setting
the constant water density $\rho = 1$, the pressure in the new nondimensional variables
is
\[
P = P_0 - gy + gp,
\]
with the nondimensional pressure variable $p$ measuring the deviation from the
hydrostatic pressure distribution. We obtain the following boundary value problem
in nondimensional variables
\begin{equation}\label{nondim}
\left\{
\begin{array}{rllllllll}
u_t+uu_x+vu_y &=& -p_x, \\[2ex]
\frac{1}{\lambda^2}[v_t+\varepsilon(uv_x+vv_y)]&=&-p_y, \\[2ex]
u_x+v_y&=& 0, \\[2ex]
p&=&\eta &\quad  \text{on} \quad y=\varepsilon \eta,\\[2ex]
v&=& \eta_t+u\eta_x &\quad  \text{on} \quad y=\varepsilon \eta,\\[2ex]
(u,v)&\to&(0,0) &\quad  \text{as} \quad y \to -\infty.
\end{array}
\right.
\end{equation}

The linearized problem is now obtained by letting $\varepsilon \to 0$ in \eqref{nondim}. 
In the limit we arrive at
\begin{equation}\label{lini}
\left\{
\begin{array}{rllllllll}
u_t &=& -p_x, \\[2ex]
\frac{1}{\lambda^2} v_t&=&-p_y, \\[2ex]
u_x+v_y&=& 0, \\[2ex]
p&=&\eta &\quad  \text{on} \quad y=0,\\[2ex]
v&=& \eta_t&\quad  \text{on} \quad y=0,\\[2ex]
(u,v)&\to&(0,0) &\quad  \text{as} \quad y \to -\infty.
\end{array}
\right.
\end{equation}
It is worth noticing that the regularity results on the streamlines for water waves obtained in the recent papers \cite{Mati} and \cite{CE3} ensure that good approximate properties hold.

Looking for solutions of \eqref{lini} representing waves traveling at speed $c > 0$, we
impose that all functions $u, v , p$ and $\eta$ have an $(x, t)$-dependence in the form of
$x-ct$. Furthermore, one seeks even periodic functions of period one, the evenness expressing the symmetry of the surface wave: it is known, cf.  \cite{CE} and  \cite{CE2, CEW} that a gravity surface traveling wave whose profile is monotone between the crests and troughs has to be symmetric about the crest.
Thus, we are
led to the fundamental Fourier mode {\it Ansatz} 
\[
\eta(x, t) = \cos [2\pi(x - ct)]. 
\]
For this  specific $\eta$ the problem \eqref{lini} has the solution 
\begin{equation*}
\left\{
\begin{array}{rllllllll}
 \eta(t,x)&=& \cos [2\pi(x - ct)]\\[1.5ex]
 u(t,x,y)&=& \frac{f(y)}{\lambda}\cos [2\pi(x - ct)],\\[1.5ex]
 v(t,x,y)&=& f(y)\sin [2\pi(x - ct)],\\[1.5ex]
 p(t,x,y)&=&  \displaystyle\frac{c}{\lambda}  \cos [2\pi(x - ct)],
\end{array}
\right.
\end{equation*}
whereby the constant $c$ and the function $f$ are given by
 $$c= \sqrt{\frac{\lambda}{2\pi}}, \quad f(y)=2\pi c\exp\left(\frac{2\pi}{\lambda}y\right).$$

Marking the variables of this nondimensional solution by tildes, we recover the physical variables by the backward
transformation
\begin{align*}
&x=\lambda\tilde {x}, \quad y= \tilde {y}, \quad t=\frac{\lambda}{\sqrt{g}}\tilde {t}, \\[1.5ex]
&u=\varepsilon\sqrt{g} \tilde{u}, \quad v= \varepsilon v\frac{\sqrt{g}}{\lambda}\tilde{v}, \quad \eta= \varepsilon\tilde{\eta},\quad p=\varepsilon \tilde{p}.
\end{align*}
Defining the wavenumber $k$ and the frequency $\omega$ by
\[
k=\frac{2\pi}{\lambda}, \quad \omega=\sqrt{gk},
\]
where $2\pi (\tilde x-c\tilde t)=kx-\omega t,$ we finally obtain that the tupel
\begin{equation}\label{linisol}
\left\{
\begin{array}{rllllllll}
\eta(t,x)&=& \varepsilon\cos(kx-\omega t),\\[1.5ex]
u(t,x,y)&=& \varepsilon \omega  \exp(ky)  \cos(kx-\omega t),\\[1.5ex]
v(t,x,y)&=& \varepsilon \omega  \exp(ky)  \sin(kx-\omega t),\\[1.5ex]
P(t,x,y)&=& P_0-gy+\varepsilon g \exp(ky)  \cos(kx-\omega t),
\end{array}
\right.
\end{equation}
is a solution of the linearised deep-water wave problem.
Note that in the physical variables the speed of the wave in \eqref{linisol} is
\[
\frac{\omega}{k}= \sqrt{\frac{g\lambda}{2\pi}}.
\]

\section{Particle trajectories}
In this last section we study the trajectories of the water particles in the linear  wave \eqref{linisol}. 
First of all, we note that a fluid particle $(x(t), y(t))$  must satisfy the following equations
\[
\frac{dx}{dt}=u, \quad \frac{dy}{dt}=v,
\]
so that, in view of \eqref{linisol}, the motion of the particle is described by the system
\begin{equation}\label{sis}
\left\{
\begin{array}{rllllllll}
\displaystyle\frac{dx}{dt}&=& M \exp(ky) \cos (kx-\omega t),\\[1ex]
\displaystyle\frac{dy}{dt}&=& M \exp(ky) \sin (kx-\omega t),\\[1ex]
\end{array}
\right.
\end{equation}
with initial data $(x(0),y(0))=(x_0, y_0).$ 
We use herein the shorthand
\begin{align}\label{em}
M:=\varepsilon \omega.
\end{align}
The right-hand side of the differential system \eqref{sis} is smooth   so that the existence
of a unique local smooth solution is ensured the Picard-Lindel\"of theorem.
Without actually solving \eqref{sis} we would like
to display the principal features of the solution. 
To study the exact solutions to \eqref{sis} it is convenient to re-write the system in a moving,
and we define therefore the variables
\begin{align}\label{transfo}
X=kx-\omega t, \quad Y=ky.
\end{align}
With this notation \eqref{sis} is equivalent to the following problem
\begin{equation}\label{smps}
\left\{
\begin{array}{rllllllll}
\displaystyle\frac{dX}{dt}&=& kM \exp(Y) \cos (X)-\omega,\\[2ex]
\displaystyle\frac{dY}{dt}&=& kM \exp(Y) \sin (X),
\end{array}
\right.
\end{equation}
with $(X(0), Y(0))=k(x_0, y_0).$
In order to determine the phase portrait corresponding to \eqref{smps}, we take into consideration  the property of \eqref{smps} of being a Hamiltonian system.
Indeed, setting
\begin{equation}\label{eq:H}
H(X,Y)=kM\exp(Y)\cos(X)-\omega Y,\quad (X,Y)\in \R^2,
\end{equation}
we re-express problem \eqref{smps} as follows
\begin{equation}\label{eq:HM}
\left\{
\begin{array}{lrrr}
\displaystyle\frac{dX}{dt}&=&\p_Y H,\\[2ex]
\displaystyle\frac{dY}{dt}&=&-\p_X H.
\end{array}
\right.
\end{equation}
Since $H$ is constant on solutions of \eqref{eq:HM}, we are left
 to determine the level curves of the function $H$.
In \cite{CEV} the authors studied the  phase portrait of \eqref{eq:HM} approximately, 
by using Morse's lemma.
A precise analysis allows us  to  determine herein  the phase portrait of the Hamiltonian system of equations describing
the motion of the particles in linear deep-water waves by using only elementary methods.

Note that $H$ is $2\pi-$periodic in $X$  and even with respect to the same variable. 
Therefore, the level curves of $H$ are determined by the level sets of the restriction  $H:\0\to \R,$
where $\0:=[0,\pi]\times\R$.
Furthermore, it is not difficult to see that the equation $\nabla H=0$ has a unique solution 
\[
P:=(X_*,Y_*):=\left(0,\ln \left( \frac{\omega}{kM}\right)\right),
\]
meaning that $P$  is the unique stationary solution of \eqref{eq:HM}  within $\0.$ 
Moreover, $P$ is located on the level curve $ H^{-1}(\{\alpha_*\}) $ with $\alpha_*$ given by
\[
\alpha_*:=\omega\left(1- \ln \left( \frac{\omega}{kM}\right) \right).
 \]

We make now a crucial observation. 
Let $\alpha\in\R$ be given with 
$H(X,Y)=\alpha$ for some $(X,Y)\in\0.$
In view of 
$kM\exp(Y)\cos(X)-\omega Y=\alpha,$
it must hold that
\[
X=f_\alpha(Y):=\arccos\left(\frac{\alpha+\omega Y}{kM\exp(Y)}\right),
\]
with $f_\alpha:D(f_\alpha)\subset\R\to[-1,1].$
The definition domain $D(f_\alpha)$ of $f_\alpha$ consists of all the points $Y\in\R$ which satisfy the inequalities
\begin{equation}\label{eq:1}
-1\leq g_\alpha(Y)=:\frac{\alpha+\omega Y}{kM\exp(Y)}\leq1.
\end{equation}
Hence $H^{-1}(\{\alpha\})\cap \0=\gr(f_\alpha) $ for all $\alpha\in\R,$
and we are left to   determine the definition domain of the function $f_\alpha $, that is to solve the inequality \eqref{eq:1}.
The first important result of this paper is the following theorem which gives a precise description of the phase curves of \eqref{eq:HM}:

\begin{figure}
$$\includegraphics[width=0.81\linewidth]{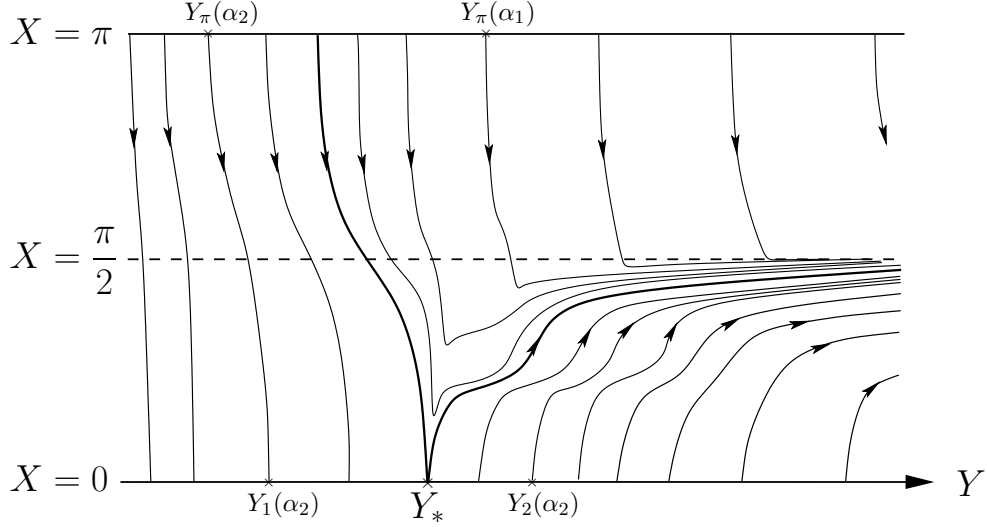}$$
\caption{Phase portrait in the moving frame ($\alpha_1<\alpha_*<\alpha_2$)}
\end{figure}

\begin{thm}\label{T}
Given $\alpha\in\R$,
there exists a  function $f_\alpha:D(f_\alpha)\subset\R\to\R$ such that
$
H^{-1}(\{\alpha\})=\gr(f_\alpha) $ for all $\alpha\in \R.$
More precisely,
there exists a unique $ Y_\pi(\alpha)>0$ such that
\begin{itemize}
\item[$(a)$] If $\alpha\leq\alpha_*$, then $D(f_\alpha)=[ Y_\pi(\alpha),\infty),$ 
\[
f_\alpha( Y_\pi(\alpha))=\pi,\quad\text{and}\quad \lim_{Y\to\infty}f_\alpha(Y)=\pi/2.
\]
Moreover, the minimum of $f_\alpha$  belongs to the interval $(0,\pi/2).$
\item[$(b)$] If $\alpha>\alpha_*$, there exist  positive real numbers $Y_i(\alpha),$ $i=1,2,$
with 
$ Y_\pi(\alpha)<Y_1(\alpha)<Y_2(\alpha)$ such that $D(f_\alpha)=[ Y_\pi(\alpha),Y_1(\alpha)]\cup[Y_2(\alpha),\infty).$
\end{itemize}
Moreover,  $f_\alpha:[ Y_\pi(\alpha),Y_1(\alpha)]\to[0,\pi]$ is bijective and decreasing, whereas $f_\alpha:[Y_2(\alpha),\infty)\to[0,\pi/2)$ is bijective and increasing.

Additionally, we have
\begin{itemize}
\item[$(1)$] $ Y_\pi:\R\to\R$ is  bijective and decreasing; 
\item[$(2)$] $ Y_1:(\alpha_*,\infty)\to(-\infty,Y_*)$ is  bijective and decreasing;
\item[$(3)$] $ Y_2:(\alpha_*,\infty)\to(Y_*,\infty)$ is  bijective and increasing;
\item[$(4)$] $ (Y_1(\alpha)-Y_\pi(\alpha))\searrow 0$ as ${\alpha\nearrow \infty}.$
\end{itemize}
The functions mentioned above are  real analytic in the interior of their definition domain (excepting $f_{\alpha_*}$ which is real analytic only on $(Y_\pi(\alpha_*), Y_*)\cup (Y_*,\infty)$).
\end{thm}

\begin{rem} 
Having proved Theorem \ref{T}, the phase portrait is complete (see Figure 2).
Note that by $(4)$, we know that the phase curves parametrised by $f_\alpha:[ Y_\pi(\alpha),Y_1(\alpha)]\to[0,\pi]$ flatten out as $\alpha\to\infty.$ 
The motion along the phase curves is determined by the fact that
\begin{align*}
-\p_X H(X,Y)=kM\exp(Y)\sin(X)>0
\end{align*}
provided $X\in(0,\pi).$
The phase portrait also discloses that the stationary solution $(X_*,Y_*)$
is  a saddle node for \eqref{eq:H}.
Solutions starting in a point $(X_0,Y_0)$ located on a separatrix (one of the level curves containing $P$)    are attracted by $(X_*,Y_*)$ provided $Y_0<Y_*.$
On the other hand, if $Y_0>Y_*$ they tend to infinity. 
\end{rem}

\begin{proof}[Proof of Theorem \ref{T}]
Let $\alpha \in \R.$ 
We observe first  that the equation  $g_\alpha(Y)=-1$ is always solvable in $\R.$
Indeed, it is equivalent to the following relation 
\begin{align*}
\alpha=h(Y):=-Mk \exp(Y)-\omega Y.
\end{align*}
Clearly, the function $h:\R\to\R$  is strictly decreasing and, moreover,
\begin{equation}\label{alf}
\lim_{Y\to-\infty}h(y)=\infty \quad\text{and} \quad  \lim_{Y\to\infty}h(y)=-\infty.
\end{equation}
Whence, there exists a unique real number  $Y_\pi(\alpha)$ such that $g_\alpha(Y_\pi(\alpha))=-1.$ 
Consequently, all phase curves of \eqref{eq:HM} intersect the line $[X=\pi]$ exactly once, that is at
\[
H^{-1}(\{\alpha\})\cap [X=\pi]=Y_\pi(\alpha) \quad \text{for all} \quad \alpha \in \R.
\]
A simple calculation shows that
\[
g'_\alpha(Y_\pi(\alpha))=\frac{1}{kM}\left(\frac{\omega}{\exp(Y_\pi(\alpha))}-\frac{\alpha+\omega Y_\pi(\alpha)}{\exp(Y_\pi(\alpha))}\right)=\frac{w+KM\exp(Y_\pi(\alpha))}{kM\exp(Y_\pi(\alpha))}>0,
\]
and we infer from the implicit function theorem that $ Y_\pi$ is real analytic and, additionally
\begin{equation}\label{eq:MM}
\frac{1}{kM\exp(Y_\pi(\alpha))}+g'_\alpha(Y_\pi(\alpha))Y_\pi'(\alpha)=0.
\end{equation}
This implies that $Y_\pi$ is strictly decreasing. 
Due to \eqref{alf}, we conclude that $Y_\pi$ is bijective.
Differentiating $g_\alpha$ at $Y\in\R,$ we observe that
\begin{equation}\label{eq:der}
g_\alpha'(Y)=\frac{\omega-\alpha-\omega Y}{kM\exp(Y)} =
\left\{ 
\begin{array}{cccc}
&\leq 0,&  \text{if}& Y\geq (\omega- \alpha)/\omega,\\[1ex]
&\geq 0, & \text{if}& Y\leq (\omega- \alpha)/\omega.
\end{array}
\right.
\end{equation} 
Let us now compare  the value of $g_\alpha$ at the threshold value $(\omega- \alpha)/\omega$, where the monotonicity  changes, to $1$, the 
larges value $g_\alpha$ can take in order to remain in the definition domain of $\arccos$.
We recover $\alpha_*$, since
\[
g_\alpha\left((\omega-\alpha)/\omega\right)\leq 1 \quad \text{if and only if}
\quad
\alpha\leq\alpha_*.
\]
Analogously, we have that $g_\alpha\left((\omega-\alpha)/\omega\right)\geq 1 $ if and only if $\alpha\geq\alpha_*.$
Since $g_\alpha((\omega-\alpha)/\omega)$ is the (positive) maximum of $g_\alpha$ and $\lim_{Y\to-\infty}g_\alpha(Y)=-\infty$, we conclude that 
\[
Y_\pi(\alpha)< (\omega-\alpha)/\omega
\]
for all $\alpha\in\R.$
Particularly, if $\alpha\leq \alpha_*$ then $g_\alpha(Y)\leq 1 $ for all $Y\geq Y_\pi(\alpha).$
Moreover, taking into account that $\lim_{Y\to\infty}g_\alpha(Y)=0,$ we conclude that $D(f_\alpha)=[Y_\pi(\alpha),\infty)$ for all $\alpha\leq\alpha_*.$ 

On the other hand, if $\alpha>\alpha_*,$ then $g_\alpha((\omega-\alpha)/\omega)>1, $  and there exist precisely two points $Y_1(\alpha)<(\omega-\alpha)/\omega<Y_2(\alpha)$ with the property that
$g_\alpha(Y_i(\alpha))=1$, $i=1,2$. 
By \eqref{eq:der} we have that $g_\alpha'(Y_1(\alpha))>0$, $g_\alpha'(Y_2(\alpha))<0,$ and we infer from the implicit function theorem 
that $Y_1$ is strictly decreasing, while $Y_2$ is strictly increasing.
We have use the fact that both $Y_1$ and $Y_2$ verify, due to $g_\alpha( Y_i(\alpha))=1$, $i=1,2,$ relation \eqref{eq:MM}. 
With this, the assertions $(1)$ and $(2)$ are immediate.

We finish the proof, by noticing that for $\alpha$ sufficiently large 
\[
2=g_\alpha(Y_1(\alpha))-g_\alpha(Y_\pi(\alpha))=g_\alpha'(Y_c(\alpha))(Y_1(\alpha)-Y_\pi(\alpha)).
\]
for some $Y_c(\alpha)\in (Y_\pi(\alpha),Y_1(\alpha)).$
Consequently, we have 
 \[
2\geq  \frac{\omega-\alpha-\omega Y_1(\alpha)}{kM\exp(Y_1(\alpha))}(Y_1(\alpha)-Y_\pi(\alpha))
\]
and, for $\alpha\to\infty,$ we obtain in view of Theorem \ref{T} $(3)$ that $Y_1(\alpha)-Y_\pi(\alpha)\to_{\alpha\to\infty}0.$
We show now that the convergence is monotonic.
Indeed, from \eqref{eq:MM} and \eqref{eq:der} we see that
\begin{align*}
Y_\pi'(\alpha)=-\frac{1}{\omega-\alpha-\omega Y_\pi(\alpha)} \qquad \text{and}\qquad Y_1'(\alpha)=-\frac{1}{\omega-\alpha-\omega Y_1(\alpha)}.
\end{align*}
 Consequently, 
 \begin{align*}
Y_1'(\alpha)-Y_\pi'(\alpha)&=\frac{\omega(Y_\pi(\alpha)-Y_1(\alpha))}{\left(\omega-\alpha-\omega Y_\pi(\alpha)\right)\left(\omega-\alpha-\omega Y_1(\alpha)\right)}\\[1ex]
&=
-\frac{2\omega}{\left(\omega-\alpha-\omega Y_\pi(\alpha)\right)\left(\omega-\alpha-\omega Y_1(\alpha)\right)}<0.
\end{align*}
This completes the proof.
\end{proof}

Knowing the phase curves $(X(t), Y(t)),$
the particle trajectories $(x(t),y(t))$ in the linear wave \eqref{linisol}
are given, via \eqref{transfo}, by 
\begin{align}\label{xy}
x(t)=\frac{X(t)}{k}+\frac{\omega}{k}t, \, \quad y(t)=\frac{Y(t)}{k}.
\end{align} 
For physical reasons realistic particle paths are bounded. 
Thus, the critical point $P$ must at least lie above the typical height of the wave.
 Since the mean surface level is zero and the amplitude is $\varepsilon,$ a crest has $y-$ coordinate  $\varepsilon.$ 
Since the $y-$ coordinate of $P$ is $\ln(1/k\varepsilon )/k$, we must thus impose that $\varepsilon k\leq \ln(1/\varepsilon k),$ that is
\begin{equation}\label{eq:cond}
 \varepsilon k \exp{(\varepsilon k)}\leq 1.
\end{equation}
This is the only condition which has to be imposed for \eqref{linisol} to be a realistic model of water waves of infinite depth, wave number $k$ and small amplitude $\varepsilon$.
This is an improvement compared to \cite{CEV}, because we do not have to to make sure that the separatrix intersects the the line $X=\pi$ only once; this is in our case a conclusion of studying the phase portrait of \eqref{eq:HM}.

We describe now the trajectory of the water particle $(x,y)$ in the fluid. 
We may consider the orbit of a particle $(x_0,y_0) $ which starts at the lowest  point (corresponding to the line $X=\pi$ in the phase-portrait).
This trajectory corresponds to the solution $(X,Y)$ of \eqref{eq:HM} which starts at $(X_0,Y_0)=(kx_0,ky_0).$
Furthermore, we denote by $\theta=\theta(Y_0)$  the time needed for $(X,Y)$ to reach the line $X=-\pi.$
Indeed, from thre first equation of \eqref{eq:HM} we have for $Y$ located beneath $Y_\pi(\alpha_*),$ that is under the separatrix, that
\begin{equation}\label{eq:theta}
\theta(Y)=2\int_0^\pi\frac{ds}{\omega-kM\exp(Y)\cos(s)}, 
\end{equation}
hence $\theta$ depends only of $Y$.
Let us observe the evolution of the particle $(x_0,y_0) $ as $t\in[0,\theta].$
Since $X=kx-\omega t,$ $Y=ky,$ we find see from \eqref{eq:HM} that

\begin{align*}
& \frac{dx}{dt}<0, \, \, \, \frac{dy}{dt}<0 \quad \text{for} \, \, X(t)\in (-\pi, -\pi/2),\\[1ex]
& \frac{dx}{dt}>0, \, \, \, \frac{dy}{dt}<0 \quad \text{for} \, \, X(t)\in (-\pi/2,0),\\[1ex]
& \frac{dx}{dt}>0, \, \, \, \frac{dy}{dt}>0 \quad \text{for} \, \, X(t)\in (0,\pi/2),\\[1ex]
& \frac{dx}{dt}<0, \, \, \, \frac{dy}{dt}>0 \quad \text{for} \, \, X(t)\in (\pi/2,\pi).
\end{align*}
Also, from the phase portrait of \eqref{eq:HM},  
we see that the hight $y$ of the particle reaches its initial value for the first time  when $t=\theta.$
What we still have to establish is the sign of $x(\theta)-x(0).$

\begin{figure}
$$\includegraphics[width=0.61\linewidth]{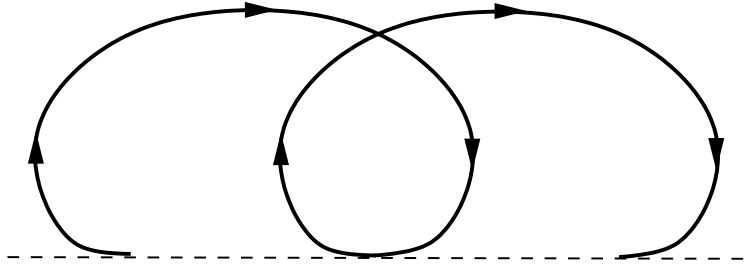}$$
\caption{Particle trajectory above the flat bed}
\end{figure}

Since
\[
x(\theta)-x(0)=\frac{\theta\omega-2\pi}{k}>0,
\]
we find from \eqref{eq:theta} that
\[
x(\theta)-x(0)=\frac{2}{k}\int_0^\pi\frac{kM\exp(Y)\cos(s)}{\omega-kM\exp(Y)\cos(s)}\, ds=\frac{2}{k}\int_0^\pi\frac{\cos(s)}{a-\cos(s)}\, ds
=\frac{4}{k}\int_0^{\pi/2}\frac{sin^2(x)}{a^2-\sin^2(x)},
\]
with $a:=\omega/(kM\exp(Y)).$
In view of $X'\leq 0,$ we conclude that the water particle experiences a forward drift equal to $x(\theta)-x(0)>0.$

Clearly, by the definition of $a$ we see that the drift decreases with depth and, in the limit $Y\to-\infty$, 
the forward drift vanishes $x(\theta)-x(0)\to0.$

Whence, we have shown that:
\begin{thm}\label{TT} The particle paths in the linear water wave \eqref{linisol} are not closed. 
There exists a forward drift over a period, and this drift decreases with depth.
\end{thm}

\end{document}